\documentclass[12pt]{article}
\usepackage{amssymb, amsmath}

\newtheorem{theorem}{Theorem}[section]

\newcommand{\vare}{\varepsilon}
\newcommand{\n}{\nonumber}

\newcommand{\si}{\sigma_R }

\newcommand{\s}{\sigma}

\renewcommand{\a}{\alpha}
\renewcommand{\o}{\omega}

\newcommand{\bb}{\begin{equation}}
\newcommand{\ee}{\end{equation}}
\newcommand{\bq}{\begin{eqnarray}}
\newcommand{\eq}{\end{eqnarray}}
\newcommand{\bqn}{\begin{eqnarray*}}
\newcommand{\eqn}{\end{eqnarray*}}

\begin{document}
\title{ Liouville-type theorems for the forced Euler equations and the Navier-Stokes equations}
\author{Dongho Chae\\
Department of Mathematics\\
            Chung-Ang University\\
          Seoul 156-756, Korea\\
              {\it e-mail : dchae@cau.ac.kr}}
 \date{}
\maketitle
\begin{abstract}
In this paper we study the Liouville-type properties for solutions
to the steady incompressible Euler equations with forces in $\Bbb
R^N$.  If we assume ``single signedness condition" on the force,
then  we can show that a $C^1 (\Bbb R^N)$ solution $(v,p)$ with
$|v|^2+ |p|\in L^{\frac{q}{2}}(\Bbb R^N)$, $q\in (\frac{3N}{N-1},
\infty)$ is trivial, $v=0$. For the solution of of the steady
Navier-Stokes equations, satisfying $v(x)\to 0$ as $|x|\to \infty$,
the condition $\int_{\Bbb R^3} |\Delta v|^{\frac65} dx<\infty$,
which is stronger than the important D-condition, $\int_{\Bbb R^3}
|\nabla v|^2 dx <\infty$, but both having the same scaling property,
implies that $v=0$. In the appendix we reprove the Theorem 1.1(\cite{cha0}), using the self-similar Euler equations directly. \\
\ \\
\noindent{\bf AMS Subject Classification Number:} 35Q30, 35Q35,
76Dxx\\
  \noindent{\bf
keywords:}   Euler equations with perturbation, steady solutions,
vanishing property
\end{abstract}
\section{ Main theorems}
 \setcounter{equation}{0}
\subsection{ The steady Euler equations with force}
Here we are concerned   on the  steady equations on $\Bbb R^N$ with
force.
  \bb\label{main}
\left\{ \aligned & ( v\cdot \nabla )  v = -\nabla p +\Phi ,\\
&\mathrm{div} \, v=0,
 \endaligned
 \right.
 \ee
 where $v=v(x)=(v_1(x),\cdots, v_N(x))$ is the velocity, and  $p=p(x)$ is the pressure.
  The force function  $\Phi[v]:\Bbb R^N\to \Bbb R^N$  satisfies the singlesignedness condition described below.
  We study Liouville-type property of the solutions to
 (\ref{main}) under this condition.
Let us fix $N\geq 2$, $k\geq 0$. Here we assume that the continuous function
$$\Phi[v] (x):=\Phi_k\left(x, v(x), Dv(x), \cdots, D^k v(x)\right)$$ for some
 $\Phi_k: \Bbb R^M \to \Bbb R^N$ for the appropriate $M(N,k)$,
 satisfies the condition of single
signedness:
 \bb\label{da1}
\mbox{either $ \Phi[v] (x)\cdot v(x)\geq 0$ or $ \Phi[v](x)\cdot
v(x)\leq 0$ for all $x\in \Bbb R^N$,}
 \ee
   and
 \bb\label{da2}
 \mbox{ $\Phi[v](x)\cdot v(x)=0$ if only if $v(x)=0$}.
 \ee
   For such
given $\Phi$ we consider the system (\ref{main}).
 Note that when $\Phi[v]=-v$ the system (\ref{main})-(\ref{da2}) becomes the usual steady
 Euler equations with a damping term. We remark that the damped
 Euler equations corresponds to a special case of the self-similar
 Euler equations(see Appendix below for more details).
 below. More generally $\Phi[v](x)= G(x, v(x),\cdots , D^k v(x))v(x)$ with a
scalar function
 $G(x, v(x),\cdots , D^k v(x))\lessgtr 0)$  satisfies (\ref{da1})-(\ref{da2}). We will prove
that a Liouville type property for the system
(\ref{main})-(\ref{da2}) under quite mild decay conditions at
infinity on the solutions. More specifically we will prove the
following.
\begin{theorem}
 Let $k\geq 0$, and $v$ be a $C^k(\Bbb R^N)$ solution of (\ref{main})-(\ref{da2}) with $\Phi=\Phi[v]$.   Suppose  there exists $q\in
(\frac{3N}{N-1} , \infty)$ such that
 \bb\label{13a}
   |v|^2+|p|\in L^{\frac{q}{2}}(\Bbb R^N).
  \ee
  Then, $v=0$.
  \end{theorem}
{\em Remark 1.1 } If $\Phi$ satisfies an extra condition
div $\Phi =0$, then the condition $p\in L^{\frac{q}{2}} (\Bbb R^N)$ can be replaced by
the well-known velocity-pressure relation in the incompressible Euler and the Navier-stokes equations,
 $$p(x)=\sum_{j,k=1}^N R_j R_k (v_jv_k)(x) $$
    with the Riesz transform $R_j$, $j=1,\cdots,N,$ in $\Bbb R^N$
    (\cite{ste}), which holds under the condition that $p(x)\to 0$ as $|x|\to \infty$.
In this case the $L^{\frac{q}{2}}$ estimate of the pressure follows from the
$L^q$ estimate for the velocity by the Calderon-Zygmund inequality,
 \bb\label{cz}\|p\|_{L^{\frac{q}{2}}}\leq C\sum_{j,k=1}^N  \|R_jR_k v_jv_k \|_{L^{\frac{q}{2}}}\leq C
 \|v\|_{L^{q}}^2 \quad 2<q<\infty.
 \ee
{\em Remark 1.2 } The theorem implies that $\mathrm{curl}
(\Phi[0]))=0$ is a necessary condition for the well-posedness of the
problem, namely $v=0$ is the unique solution of the equations.

\subsection{ The steady Navier-Stokes equations in $\Bbb R^3$}
  Here we study the following  system of steady Navier-Stokes equations in $\Bbb R^3$.
 $$
(NS)\left\{ \aligned & ( v\cdot \nabla )  v = -\nabla p +\Delta v,\\
&\mathrm{div} \, v=0,
 \endaligned
 \right.
 $$
We consider here the generalized solutions of the system (NS),
satisfying
 \bb\label{diri}
\int_{\Bbb R^3} |\nabla v|^2dx <\infty,
 \ee
 and
    \bb\label{12c}
  \lim_{|x|\to \infty} v(x)= 0.
  \ee
 It is well-known that a generalized solution to (NS) belonging to $ W^{1,2}_{loc}(\Bbb R^3)$ implies that
 $v$ is smooth(see e.g.\cite{gal}). Therefore without loss of
 generality we can assume that our solutions to (NS) satisfying (\ref{diri}) are smooth.
 The uniqueness question, or equivalently the question of  Liouville
 property of solution for the system (NS) under the assumptions
 (\ref{diri}) and (\ref{12c}) is a long standing open problem.
  On the other hand, it is well-known that the uniqueness of solution holds in the
   class $L^{\frac92} (\Bbb R^3)$, namely a
  smooth solution to (NS) satisfying
$v\in L^{\frac92} (\Bbb R^3)$ is $v=0$(see Theorem 9.7 of \cite{gal}). We assume here slightly stronger condition than
  (\ref{diri}), but having the same scaling property, to deduce our Liouville-type result.
\begin{theorem}
Let $v$ be a smooth solution of (NS) satisfying (\ref{12c}) and
 \bb\label{sdiri}
 \int_{\Bbb R^3} |\Delta v|^{\frac65}\, dx<\infty.
 \ee
   Then, $v=0$ on
 $\Bbb R^3$.
  \end{theorem}
  {\em Remark 1.3 } Under the assumption (\ref{12c}) we have the
  inequalities with the norms of the {\em same scaling properties,}
$$
 \|v\|_{L^6}\leq C\|\nabla v\|_{L^2} \leq C  \|D^2 v\|_{L^{\frac65}}
 \leq C\|\Delta v\|_{L^\frac65}<\infty
$$
due to the Sobolev and the Calderon-Zygmund inequalities. Thus,
(\ref{sdiri}) implies (\ref{diri}).  There is no, however, mutual
implication relation between Theorem 1.2 and
the above mentioned $L^{\frac92}$ result, although our assumption (\ref{sdiri}) corresponds to $L^6(\Bbb R^3)$
at the level of scaling. \\

\section{Proof of the Main Theorems }
\setcounter{equation}{0}

\noindent{\bf Proof of Theorem 1.1 }  We denote
 $$[f]_+=\max\{0, f\}, \quad  [f]_-=\max\{0, -f\},$$
and
$$
D_\pm:=\left\{ x\in \Bbb R^N\, \Big|\,\left[p(x)+\frac12
|v(x)|^2\right]_\pm >0\right\}
$$
respectively. We introduce the radial cut-off function $\sigma\in
C_0 ^\infty(\Bbb R^N)$ such that
 \bb\label{16}
   \sigma(|x|)=\left\{ \aligned
                  &1 \quad\mbox{if $|x|<1$},\\
                     &0 \quad\mbox{if $|x|>2$},
                      \endaligned \right.
 \ee
and $0\leq \sigma  (x)\leq 1$ for $1<|x|<2$.  Then, for each $R
>0$, we define
 $$
\s \left(\frac{|x|}{R}\right):=\s_R (|x|)\in C_0 ^\infty (\Bbb R^N).
$$
  We multiply first equations of (\ref{main}) by $v$ to obtain
   \bb\label{main1}
   v\cdot \Phi =v\cdot \nabla \left(p+\frac12 |v|^2\right).
   \ee
   Next, we multiply (\ref{main1}) by
  $ \left[p+\frac12 |v|^2\right]_+
^{\frac{qN-q-3N}{2N}}\s_R \, \mathrm{sign}\{v\cdot \Phi \} $ and
integrate over $\Bbb R^N$ to have
  \bq\label{euler1}
 \lefteqn{\int_{\Bbb R^N} \left[p+\frac12
|v|^2\right]_+^{\frac{qN-q-3N}{2N}}\left|v\cdot \Phi\right| \s_{R}\,
dx}\n \\
&&=\mathrm{sign}\{ v\cdot \Phi\}\int_{\Bbb R^N} \left[p+\frac12
|v|^2\right]_+^{\frac{qN-q-3N}{2N}}\s_R v \cdot\nabla
\left(p+\frac12 |v|^2\right) \,dx\n \\
&&:=I
   \eq
   We estimate $I$ as follows.
 \bqn
 |I|&=&\left|\int_{\Bbb R^N} \left[p+\frac12 |v|^2\right]_+^{\frac{qN-q-3N}{2N}} \s_R
v\cdot\nabla \left(p+\frac12 |v|^2 \right)\, dx\right| \n \\
&=&\left|\int_{D_+} \left[p+\frac12
|v|^2\right]_+^{\frac{qN-q-3N}{2N}} \s_R
v\cdot\nabla \left[p+\frac12 |v|^2 \right]_+\, dx\right|\n \\
 &=& \frac{2N}{qN-q-N} \left|\int_{D_+}\si v\cdot\nabla \left[p+\frac12
|v|^2 \right]_+^{\frac{qN-q-N}{2N}} \, dx \right|\n \\
&=&\frac{2N}{qN-q-N} \left| \int_{D_+}\left[p+\frac12 |v|^2
\right]_+^{\frac{qN-q-N}{2N}}  v\cdot\nabla \si \, dx\right| \n \\
 &\leq&
\frac{C\|\nabla\s\|_{L^\infty} }{R} \left(\int_{\Bbb R^N}
(|p|+|v|^2)^{\frac{q}{2}} \, dx\right)^{\frac{qN-q-N}{qN}}
\|v\|_{L^q(R\leq |x|\leq 2R)} \times\n \\
&&\hspace{.5in} \times\left(\int_{\{ R\leq |x|\leq 2R\}} \,
dx \right)^{\frac1N}\n \\
  &\leq& C\|\nabla\s\|_{L^\infty}
\left(\|p\|_{L^{\frac{q}{2}}} +\|v\|_{L^q}^2
\right)^{\frac{qN-q-N}{qN}}\|v\|_{L^q(R\leq |x|\leq 2R)}\to 0
 \eqn
as $R\to \infty$.
 Therefore, passing $R\to \infty$ in (\ref{euler1}), we obtain
 \bb\label{116a}
 \int_{\Bbb R^N} \left[p+\frac12
|v|^2\right]_+^{\frac{qN-q-3N}{2N}} \left|v\cdot \Phi\right| \, dx
=0 \ee by the Lebesgue Monotone Convergence Theorem.
  Similarly, multiplying (\ref{main1}) by $ \left[p+\frac12
  |v|^2\right]_-
^{\frac{qN-q-3N}{2N}} \s_R $, and integrate over $\Bbb R^N$,  we
deduce by similarly to the above,
 \bq\label{116aa}
 \lefteqn{\int_{\Bbb R^N}
 \left[p+\frac12 |v|^2\right]_-^{\frac{qN-q-3N}{2N}} \left|v\cdot \Phi\right| \s_{R}\, dx}\hspace{.3in}\n \\
 &&=-\int_{\Bbb R^N} \left[p+\frac12 |v|^2\right]_-^{\frac{qN-q-3N}{2N}} \s_R
v\cdot\nabla \left(p+\frac12 |v|^2 \right)\, dx \n \\
&&=\int_{\Bbb R^N} \left[p+\frac12
|v|^2\right]_-^{\frac{qN-q-3N}{2N}} \s_R
v\cdot\nabla \left[p+\frac12 |v|^2 \right]_-\, dx\n \\
&&\leq C\|\nabla\s\|_{L^\infty}\left(\|p\|_{L^{\frac{q}{2}}}
+\|v\|_{L^q}^2 \right)^{\frac{qN-q-N}{qN}}\|v\|_{L^q(R\leq |x|\leq
2R)}\to 0\n \\
 \eq
as $R\to \infty$. Hence,
 \bb\label{116ab}\int_{\Bbb R^N}
\left[p+\frac12 |v|^2\right]_-^{\frac{qN-q-3N}{2N}}\left|v\cdot
\Phi\right|\, dx =0 \ee by the Lebesgue Monotone Convergence Theorem
again.
 Let us define
$$ \mathcal{S}=\{ x\in \Bbb R^N\, |\, v(x)\neq0\}. $$
We note that $\mathcal{S}$ is an open set in $\Bbb R^N$.
Suppose $\mathcal{S}\neq \emptyset$. Then, (\ref{116aa}) and
(\ref{116ab}) together with (\ref{da1})-(\ref{da2}) imply
$$
\left[p(x)+\frac12 |v(x)|^2\right]_+=\left[p(x)+\frac12
|v(x)|^2\right]_-=0\quad \forall x\in \mathcal{S}.
$$
Namely,
$$
p(x)+\frac12 |v(x)|^2=0\quad \forall x\in \mathcal{S}.
$$
Since this holds for any open subset of $\mathcal{S}$, we have also
$ \nabla (p+\frac12 |v|^2) (x)=0$ for all $x\in \mathcal{S}$.
From (\ref{main1}) this implies
 \bb\label{119}
  \Phi[v](x)\cdot v(x)=0\qquad \forall x\in \mathcal{S}.
  \ee
 Considering the conditions on $\Phi$ in (\ref{da1})-(\ref{da2}), we have a contradiction, and therefore we need $\mathcal{S}=\emptyset$, namely
 $v=0$ on $\Bbb R^N$.
   $\square$\\
\ \\
Next in order to prove Theorem 1.2 we recall the following result proved by Galdi(see Theorem X.5.1 of \cite{gal} for more general version).
\begin{theorem}
Let $v(x)$ be a generalized solution of (NS) satisfying (\ref{diri}) and (\ref{12c}) and $p(x)$ be the associated pressure, then there exists $p_1\in \Bbb R$ such that
$$\lim_{|x|\to \infty} |D^\alpha v(x)|+\lim_{|x|\to \infty} |D^\alpha \left(p(x)-p_1\right)|=0 $$
uniformly for all multi-index $\alpha=(\alpha_1, \alpha_2, \alpha_3)\in [\Bbb N\cup\{0\}]^3$.
\end{theorem}

\noindent{\bf Proof of Theorem 1.2 }
 Under the
assumption (\ref{sdiri}) and Remark 1.1, Theorem IX.6.1 of
\cite{gal} implies
 that
  \bb\label{decay}
 \lim_{|x|\to \infty}|p(x)-p_1 |=0.
  \ee
 for  a constant $p_1$.
 Therefore, if we set
 $$
  \label{ber}Q(x):=\frac12 |v(x)|^2 +p(x)-p_1,
 $$
 then
  \bb\label{13}
\lim_{|x|\to \infty} |Q(x)|=0.
  \ee
As before we  denote $[f]_+=\max\{0, f\}, \quad  [f]_-=\max\{0,
-f\}.$ Given $\vare > 0$, we define
  \bqn
    D_+^\vare&=&\left\{ x\in \Bbb
R^3\, \Big|\,\left[Q(x)-\vare\right]_+>0\right\},\n
\\
D_-^\vare&=&\left\{ x\in \Bbb R^3\,
\Big|\,\left[Q(x)+\vare\right]_->0\right\}.
  \eqn
 respectively. Note
that (\ref{13}) implies that $D_\pm^\vare$ are bounded sets in $\Bbb
R^3$. Moreover,
  \bb\label{obs} Q\mp\vare
=0\quad\mbox{on}\quad \partial D_\pm^\vare
 \ee
 respectively.
Also, thanks to the Sard theorem combined with the implicit function
theorem $\partial D_\pm^\vare$'s are smooth level surfaces in $\Bbb
R^3$ except the values of $\vare>0$, having the zero Lebesgue
measure, which corresponds to the critical values of $z=Q(x)$. It is
understood that our values of  $\vare$ below avoids these
exceptional ones. We write the system (NS) in the form,
 \bb\label{17}
 -v\times\mathrm{ curl} \, v =-\nabla Q +\Delta v.
 \ee
 Let us
multiply (\ref{17}) by $ v \left[Q-\vare\right]_+$, and integrate it
over $\Bbb R^3$. Then, since
  $
v\times \mathrm{curl}\,v  \cdot v =0, $
  we have
 \bq\label{18}
0&=& -\int_{\Bbb R^3} \left[Q-\vare\right]_+ v \cdot\nabla
\left(Q-\vare\right) \,dx+\int_{\Bbb R^3}\left[Q-\vare\right]_+ v\cdot \Delta v\, dx\n\\
&:=&I_1 +I_2 .
 \eq
Integrating by parts, using (\ref{obs}), we obtain
 \bqn
 I_1=-\int_{D_+^\vare} \left(Q-\vare\right)
v\cdot\nabla \left(Q -\vare\right)\, dx= -\frac{1}{2}
\int_{D_+^\vare}
  v\cdot\nabla \left(Q-\vare \right)^{2} \, dx.
=0 \eqn
 Using
 \bb\label{f2aa}
 v\cdot \Delta v=\Delta (\frac12|v|^2)-|\nabla v|^2,
 \ee
 and the well-known formula for the Navier-Stokes equations,
 \bb\label{f2a}
 \Delta p=|\o|^2-|\nabla v |^2,
 \ee
 we have
 \bq\label{110}
 I_2&=&-\int_{\Bbb R^3} |\nabla v|^2\left[ Q-\vare\right]_+ \, dx +\int_{\Bbb R^3}\Delta \left(\frac12 |v|^2\right)
\left[Q-\vare\right]_+ \, dx \n \\
  &=&-\int_{\Bbb R^3}|\o|^2\left[Q-\vare\right]_+\, dx +\int_{\Bbb R^3} \Delta \left(Q-\vare\right)
 \left[ Q-\vare\right]_+
 \, dx\n \\
&:=&J_1+J_2.
 \eq
Integrating by parts, we transform $J_2$ into
 \bq\label{110a}
 J_2=\int_{D_+^\vare} \Delta\left(Q-\vare\right)
 \left( Q-\vare\right)
 \,dx=-
 \int_{D_+^\vare}\left|\nabla\left(Q-\vare\right)\right|^2
 \,dx.
 \eq
Thus,  the derivations (\ref{18})-(\ref{110a})
 lead us to
 \bb\label{dom}
0=\int_{D_+^\vare} |\o|^2\left| Q-\vare\right|\,
dx+\int_{D_+^\vare}\left|\nabla\left(Q-\vare\right)\right|^2
 \,dx
\ee for all $\vare>0$. The vanishing of the second term of
(\ref{dom}) implies
$$
\left[Q-\vare\right]_+=C_0\quad \mbox{on}\quad D_+^\vare
$$ for a constant $C_0$. From the fact (\ref{obs})
 we have $C_0=0$, and $[Q-\vare]_+=0$ on $\Bbb R^3$, which holds for all  $\vare >0.$
Hence,
 \bb\label{plus}
\left[Q\right]_+=0 \quad\mbox{on $\Bbb R^3$}.
  \ee
This shows that $Q\leq 0$ on $\Bbb R^3$. Suppose $Q=0$ on $\Bbb
R^3$. Then, from (\ref{17}), we have $v\cdot\Delta v=0$ on $\Bbb
R^3$. Hence,
$$ \Delta p=-\frac12\Delta |v|^2=-v\cdot\Delta v-|\nabla v|^2=-|\nabla
v|^2.
$$
Comparing this with (\ref{f2a}), we have $\o=0$. Combining this with
div $v=0$, we find that $v$ is a harmonic function in $\Bbb R^3$.
Thus, by (\ref{12c}) and the Liouville theorem for the harmonic
function, we have  $v=0$, and we are done. Hence, without loss of generality, we may assume
$$0>\inf_{x\in \Bbb R^3} Q(x).
$$
Given $\delta>0$, we multiply (\ref{17}) by $ v
\left[Q+\vare\right]_- ^{\delta} $, and integrate it over $\Bbb
R^3$. Then, similarly to the above
  we have
 \bq\label{m18}
0&=& -\int_{\Bbb R^3} \left[Q+\vare\right]_-^{\delta} v \cdot\nabla
\left(Q+\vare\right) \,dx+\int_{\Bbb R^3} \left[Q+\vare\right]_-^{\delta} v\cdot \Delta v\, dx\n\\
&:=&I_1' +I_2' .
 \eq
Observing $Q(x)+\vare=- \left[Q(x)+\vare\right]_-$ for all $x\in
D_-^\vare$, integrating by part, we obtain
 \bqn
 I_1'&=&\int_{D_-^\vare}\left[Q+\vare\right]_-^{\delta}
v\cdot\nabla \left[Q+\vare\right]_-\, dx\n \\
 &=&\frac{1}{\delta+1} \int_{D_-^\vare} v\cdot\nabla \left[Q+\vare \right]_-^{\delta+1} \, dx =0.
 \eqn
 Thus, using (\ref{f2aa}), we have
 \bq\label{m110}
 0=-\int_{D_-^\vare} |\nabla v|^2\left[Q+\vare\right]_-^{\delta} \,
 dx+\frac12 \int_{D_-^\vare}\left[ Q+\vare\right]_- ^{\delta}\Delta |v|^2
\, dx \eq
 Now,  we have the point-wise convergence
 $$\left[Q(x)+\vare\right]_-^\delta \to 1 \quad \forall x\in D_-^\vare.
 $$
 as $\delta\downarrow 0$. Since
 \bqn
\int_{\Bbb R^3}| v\cdot \Delta v|\, dx&\leq& \|v\|_{L^6}\|\Delta
v\|_{L^{\frac65}}\leq C \|\nabla
v\|_{L^2}\|\Delta v\|_{L^{\frac65}}\n \\
&\leq &C \|\Delta v\|_{L^{\frac65}}^2<\infty,
 \eqn
 we have
 \bb\label{elone}\Delta |v|^2=2 v\cdot \Delta v +2 |\nabla v|^2 \in L^1(\Bbb R^2).
 \ee
 Hence, passing $\delta\downarrow 0$ in
 (\ref{m110}),
 by the dominated convergence
 theorem, we obtain
 \bb\label{m110a}
 \int_{D_-^\vare} |\nabla v|^2\,dx= \frac12\int_{D_-^\vare}
  \Delta |v|^2\, dx,
 \ee
which holds for all $\vare>0$. For a sequence $\{\vare_n\}$ with
$\vare_n \downarrow 0$ as $n\to \infty$, we observe
 $$ D_-^{\vare_n }\uparrow \cup_{n=1}^\infty D_-^{\vare_n}=D_-:=
 \{ x\in \Bbb R^3\, |\, Q(x)<0.\}.
 $$
Thus, observing (\ref{elone}) again, we can apply the dominated
convergence theorem in passing $\vare \downarrow 0$ in (\ref{m110a})
to deduce
 \bb\label{mmdom1}
\int_{D_-} |\nabla v|^2\, dx=\frac12\int_{D_-}\Delta |v|^2\, dx.
 \ee
 Now, thanks to (\ref{plus}) the set
 $$ S=\{ x\in \Bbb R^3\, |\, Q(x)=0\} $$
 consists of critical(maximum) points of $Q$, and hence
 $ \nabla Q(x)=0$ for all $x\in S,$ and the system (\ref{17}) reduces to
  \bb\label{redu}
   -v\times \o=\Delta v \quad\mbox{on}\quad S.
  \ee
Multiplying (\ref{redu}) by $v$, we have that
$$ 0=v\cdot \Delta v=\frac12\Delta |v|^2-|\nabla v|^2\quad\mbox{on}\quad S.
$$
Therefore, one can extend the domain of integration in
(\ref{mmdom1}) from $D_-$ to $D_-\cup S=\Bbb R^3$, and therefore
 \bb\label{mmdom2}
\int_{\Bbb R^3} |\nabla v|^2\, dx=\frac12\int_{\Bbb R^3}\Delta
|v|^2\, dx.
 \ee
 We now claim the right hand side of (\ref{mmdom2}) vanishes.
 Indeed, since $\Delta |v|^2\in L^1(\Bbb R^3)$ from (\ref{elone}), applying the dominated
 convergence theorem, we have
 \bqn
 \left|\int_{\Bbb R^3} \Delta |v|^2\, dx\right|&=&\lim_{R\to \infty} \left|\int_{\Bbb R^3}
 \Delta |v|^2 \si \, dx\right|
 =\lim_{R\to \infty} \left|\int_{\Bbb R^3}
  |v|^2 \Delta\si \, dx\right|\n \\
  &\leq&\lim_{R\to \infty}\int_{\Bbb R^3}
  |v|^2 |\Delta\si| \, dx\n \\
 &\leq& \lim_{R\to \infty}\frac{\|D^2\s\|_{L^\infty}}{R^2}
 \|v\|_{L^6(R\leq |x|\leq
 2R)}^2
 \left(\int_{\{R\leq |x|\leq
 2R\}} dx\right)^{\frac23}\n \\
 &\leq &C\|D^2\s\|_{L^\infty}\lim_{R\to \infty}
 \|v\|_{L^6(R\leq |x|\leq
 2R)}^2=0
 \eqn
 as claimed.
 Thus (\ref{mmdom2}) implies that
 $$
\nabla v=0 \quad\mbox{on} \quad \Bbb R^3,
$$
and $v=$ constant. By (\ref{12c}) we have $v=0$.
$\square$\\
\ \\
 \noindent{\em Remark after the proof of Theorem 1.2: } The first part of the
above proof, showing $ [Q]_+=0 $ can be also done by applying the
maximum principle, which follows from the following identity for $Q$,
$$ -\Delta Q+v\cdot \nabla Q =-|\o|^2 \leq 0
$$
I do not think, however, the maximum principle can also be applied
to the proof of the second part, showing $[Q]_-=0$, which is more
subtle than the first part. The above proof overall shows that the
argument of the proof I used for this second part can also be
adapted for the first part without using the
maximum principle, which exhibits consistency.\\

\[ \mbox{\large\bf Appendix}\]

\appendix
\section{Remarks on the self-similar Euler equations}
\setcounter{equation}{0}
\numberwithin{equation}{section}
Let $a,b$ are given constants with $b\neq 0$. We study here the system in $\Bbb R^3$.
 \bb\label{pe}
\left\{ \aligned & ( v\cdot \nabla )  v = -\nabla p  +a v+b(x\cdot \nabla) v,\\
&\mathrm{div} \, v=0.
 \endaligned
 \right.
 \ee
In the special case of $a=-\frac{\alpha}{\alpha+1}$,
$b=-\frac{1}{\alpha+1}$ the system (\ref{pe})  reduces to the
self-similar Euler equations.
 \bb\label{ss}
(SSE) \left\{ \aligned & \frac{\alpha}{\alpha+1} v+\frac{1}{\alpha+1}(x\cdot \nabla) v +( v\cdot \nabla )  v = -\nabla p ,\\
&\mathrm{div} \, v=0,
 \endaligned
 \right.
 \ee
 which is obtained from the time dependent Euler equations,
 $$
 \left\{ \aligned & u_t+(u\cdot \nabla )u =-\nabla \pi\\
&\mathrm{div} \, u=0,
 \endaligned
 \right.
$$
by the self-similar ansatz, \bqn
u(x,t) & = &\frac{1}{(T - t)^{\frac{\a}{1+\a}}} v\left( \frac{x-x^*}{(T - t)^{\frac{1}{1+\a}}}\right), \\
\pi(x,t) & = & \frac{1}{(T - t)^{\frac{2\a}{1+\a}}} p\left(
\frac{x-x^*}{(T - t)^{\frac{1}{1+\a}}}\right). \eqn Note that the
damped Euler equation, which is a trivial case of (\ref{main}) is
the case when $\alpha=\infty$ in (\ref{ss}).
 In \cite{cha0}, in particular, a Liouville-type result
 for the system (\ref{ss}) was derived, using the time dependent Euler equations, where
 we need to use existence result of a back-to-label map due to Constantin(\cite{con0}).
 In the following we prove similar result for the general system (\ref{ss}).
\begin{theorem} Let $v$ be a $C^2 (\Bbb R^3)$ solution to (\ref{pe}) with $b\neq 0$, satisfying
 \bb\label{condi13} \|\nabla v\|_{L^\infty}<\infty
\quad\mbox{and}\quad \o\in \bigcup_{r >0}\bigcap_{0<q< r} L^q(\Bbb R^3).
  \ee
 Then, $v=\nabla h$ for a  harmonic scalar function $h$
on $\Bbb R^3$. Thus, if we impose further the condition $
\lim_{|x|\to \infty} |v(x)|=0$, then $v=0$.
\end{theorem}
 \noindent{\bf Proof  }
 We first observe that from the calculus identity
 $$ v(x)=v(0)+\int_0 ^1\partial_s v(sx) ds=v(0)+\int_0 ^1 x\cdot
 \nabla v(sx) ds,
 $$
 we have
 $|v(x)|\leq |v(0)|+ |x|\|\nabla v\|_{L^\infty}\leq C(1+|x|)(\|\nabla
 v\|_{L^\infty}+|v(0)|),$
 and
  \bb\label{cal}
  \sup_{x\in \Bbb R^3} \frac{|v(x)|}{1+|x|} \leq C(\|\nabla
  v\|_{L^\infty}+|v(0)|).
  \ee
We consider the vorticity equation of (\ref{pe}),
 \bb\label{vor}
 -(a+b)\o- b( x\cdot \nabla ) \o+(v\cdot \nabla
 )\o=(\o\cdot \nabla )v .
 \ee
 Let $\delta >0$, and  take $L^2(\Bbb R^3)$ inner product (\ref{vor}) by $\o (\delta +|\o|^2)^{\frac{q}{2}-1}
 \si$, and integrate over $\Bbb R^3$ to obtain
 \bq\label{newvor}
&& -(a+b)\int_{\Bbb R^3} |\o|^2(\delta +|\o|^2)^{\frac{q}{2}-1}\si dx - \int_{\Bbb R^3} [\o\cdot \nabla )v]\cdot\o (\delta +|\o|^2)^{\frac{q}{2}-1}\si dx\n \\
&&\qquad = \frac{1}{q}\int_{\Bbb R^3} \left[(( bx-v)\cdot \nabla
)(\delta +|\o|^2)^{\frac{q}{2}}\right]\si dx.\n \\
 \eq
 For fixed $\delta>0$ and $R>0$ the integrands in the right hand side of (\ref{newvor}) are sufficiently smooth functions having the compact support, and one can integrate by part  them to obtain
 \bq\label{newvor1}
&& -(a+b)\int_{\Bbb R^3} |\o|^2(\delta +|\o|^2)^{\frac{q}{2}-1}\si dx - \int_{\Bbb R^3}[ (\o\cdot \nabla )v]\cdot\o (\delta +|\o|^2)^{\frac{q}{2}-1}\si dx\n \\
&&\qquad = -\frac{3b}{q}\int_{\Bbb R^3}(\delta +|\o|^2)^{\frac{q}{2}}\si dx
-\frac{1}{q}\int_{\Bbb R^3} (\delta +|\o|^2)^{\frac{q}{2}} \left(( bx-v)\cdot \nabla \right)\si dx.\n \\
 \eq
 Passing $\delta \downarrow 0$ in (\ref{newvor1}), using the dominated convergence theorem, we have
 \bq\label{sel}
 \lefteqn{\left(-a-b+ \frac{3b}{q}\right)\int_{\Bbb R^3} |\o|^q \si dx-
\int_{\Bbb R^3} (\o \cdot \nabla) v\cdot \o  |\o|^{q-2}\si\, dx}\n \\
&&=-\frac{ b}{q} \int_{\Bbb R^3} |\o|^q (x\cdot\nabla )\si \,
dx +\frac{1}{q} \int_{\Bbb R^3} |\o|^q (v\cdot\nabla )\si \, dx\n \\
&&:=I+J.
 \eq
We estimate $I$ and $J$ easily as follows.
  $$
  |I|\leq  \frac{|b|}{ qR}\int_{\{R\leq |x|\leq 2R\}} |\o|^q
  |x||\nabla \s|\, dx
  \leq \frac{2|b|}{q}\|\nabla \s\|_{L^\infty}\|\o\|^q_{L^p(R\leq |x|\leq
  2R)}\to 0
  $$
   as $R\to \infty$.
   \bqn
    |J|&\leq&  \frac{1}{q R}\int_{\{R\leq |x|\leq 2R\}} |\o|^q
  |v||\nabla \s|\, dx
  \leq\frac{1+2R}{q R}\int_{\{R\leq |x|\leq 2R\}}
  \frac{|v(x)|}{1+|x|}|\o|^q|\nabla \s|\, dx\n \\
  &\leq &\frac{C(1+2R)}{q R}\|\nabla \s\|_{L^\infty}(\|\nabla v\|_{L^\infty}+|v(0)|) \|\o\|^q_{L^p(R\leq |x|\leq
  2R)}\to 0
  \eqn
  as $R\to \infty$, where we used (\ref{cal}). Therefore, passing
  $R\to \infty$ in (\ref{sel}), and using the dominated convergence theorem for the left hand side, we obtain,
  $$
  \left(-a-b+ \frac{3b}{q}\right)\int_{\Bbb R^3} |\o|^q  dx=
\int_{\Bbb R^3} (\o \cdot \nabla) v\cdot \o  |\o|^{q-2}\, dx,
  $$
  from which we deduce easily
  \bb\label{sel2} -\|\nabla v\|_{L^\infty}\|\o\|_{L^q}^q
\leq \left(-a-b+\frac{3b}{q}\right)\|\o\|_{L^q}^q\leq \|\nabla
v\|_{L^\infty}\|\o\|_{L^q}^q.
  \ee
 Suppose there exists $x_0\in \Bbb R^3$ such that  $\o(x_0)\neq 0$, then since $\o$ is a continuous function, one has $\|\o\|_{L^q} > 0$, and we can divide (\ref{sel2}) by
 $\|\o\|_{L^q}^q$ to have
\bb\label{sel3}
 -\|\nabla v\|_{L^\infty} \leq
\left(-a-b+\frac{3b}{q}\right)\leq \|\nabla v\|_{L^\infty},
  \ee
 which holds for all $q\in (0, r)$ and for some $r>0$. Since $b\neq 0$, passing $q\downarrow 0$ in (\ref{sel3}), we
 obtain desired contradiction. Therefore $\o=\mathrm{curl}\, v=0$. This,
 together with $\mathrm{div}\,v=0$, provides us with the fact that $v=\nabla
 h$ for a scalar harmonic function $h$ on $\Bbb R^3$. $\square$\\

 $$ \mbox{\bf Acknowledgements} $$
 The author would like to thank deeply to the anonymous referee for careful
reading and constructive criticism. This work was supported
partially by the NRF grant. no. 2006-0093854 and also by Chung-Ang
University Research Grants in 2012.

\end{document}